\documentclass[a4paper,UKenglish,cleveref, autoref, thm-restate]{lipics-v2021}

\nolinenumbers
\usepackage{amsfonts, amsmath, amssymb, mathtools, extarrows, bbm, fancyvrb}
\usepackage[small,nohug,heads=LaTeX]{diagrams}
\diagramstyle[labelstyle=\scriptstyle]

\newcommand{\ra}{\rightarrow}

\newcommand{\mcb}{\mathcal{B}}
\newcommand{\mcc}{\mathcal{C}}

\newcommand{\mco}{\mathcal{O}}
\newcommand{\mcp}{\mathcal{P}}

\newcommand{\mbb}{\mathbb{B}}

\newcommand{\mbu}{\mathbb{U}}

\newcommand{\nats}{\mathbb{N}}

\newcommand{\stl}{\left\{}
\newcommand{\str}{\right\}}

\newcommand{\bu}{\bullet} 
\newcommand{\hm}{\text{Hom}}

\newcommand{\und}{\underline}
\newcommand{\mbbo}{\mathbbm{1}}
\newcommand{\sets}{\textbf{Sets}}

\newcommand{\el}{\textbf{El}}
\DeclareSymbolFont{stmry}{U}{stmry}{m}{n}
\DeclareMathSymbol\llparenthesis\mathop{stmry}{"4C}
\DeclareMathSymbol\rrparenthesis\mathop{stmry}{"4D}
\newcommand{\type}{{\color{teal}\textbf{Type}}}
\newcommand{\typebl}{\textbf{Type}}
\newcommand{\prop}{{\color{teal}\textbf{Prop}}}
\newcommand\blfootnote[1]{%
  \begingroup
  \renewcommand\thefootnote{}\footnote{#1}%
  \addtocounter{footnote}{-1}%
  \endgroup
}

\title{A Formalization of Operads in Coq} 

\author{Zachary Flores}{High Assurance Soluions, Two Six Technologies}{zachary.flores@twosixtech.com}{}{}

\author{Angelo Taranto}{High Assurance Solutions, Two Six Technologies}{angelo.taranto@twosixtech.com}{}{}

\author{Eric Bond}{High Assurance Solutions, Two Six Technologies}{eric.bond@twosixtech.com}{}{}

\author{Yakir Forman}{High Assurance Solutions, Two Six Technologies}{yakir.forman@twosixtech.com}{}{}

\authorrunning{J. Open Access and J.\,R. Public} 

\Copyright{Zachary Flores, Angelo Taranto, Eric Bond, and Yakir Forman} 

\ccsdesc[500]{Theory of computation~Constructive mathematics}

\keywords{Operads, Formal Mathematics, Coq} 

\supplement{Disclaimer:  This research was developed with funding from the Defense Advanced Research Projects Agency
(DARPA). The views, opinions and/or findings expressed are those of the
author and should not be interpreted as representing the official views or
policies of the Department of Defense or the U.S. Government.}

\funding{Supported by DARPA V-SPELLS HR001120S0058.}

\acknowledgements{We would like to thank Li-yao for their help in writing Coq code for multi-composition.}

\EventEditors{John Q. Open and Joan R. Access}
\EventNoEds{2}
\EventLongTitle{42nd Conference on Very Important Topics (CVIT 2016)}
\EventShortTitle{CVIT 2016}
\EventAcronym{CVIT}
\EventYear{2016}
\EventDate{December 24--27, 2016}
\EventLocation{Little Whinging, United Kingdom}
\EventLogo{}
\SeriesVolume{42}
\ArticleNo{23}

\begin{document}

\maketitle

\begin{abstract}


What provides the highest level of assurance for correctness of execution within a programming language?  One answer, and our solution in particular, to this problem is to provide a formalization for, if it exists, the denotational semantics of a programming language.  Achieving such a formalization provides a gold standard for ensuring a programming language is correct-by-construction.  In our effort on the DARPA V-SPELLS program, we worked to provide a foundation for the denotational semantics of a meta-language using a mathematical object known as an operad.  This object has compositional properties which are vital to building languages from smaller pieces. In this paper, we discuss our formalization of an operad in the proof assistant Coq.  Moreover, our definition within Coq is capable of providing proofs that objects specified within Coq are operads.  This work within Coq provides a formal mathematical basis for our meta-language development within V-SPELLS.  Our work also provides, to our knowledge, the first known formalization of operads within a proof assistant that has significant automation, as well as a model that can be replicated without knowledge of Homotopy Type Theory.    

\end{abstract} 

\section{Introduction}
\blfootnote{\textbf{Distribution Statement A}:  Approved for Public Release, Distribution Unlimited}

 The DARPA V-SPELLS (Verified Security and Performance Enhancement of Large Legacy Software) program aims to create developer-accessible capability for piece-by-piece enhancement of software components for large legacy codebases with new verified code that is \textit{safely composable }with the rest of the system.  

In our approach with the Johns Hopkins Applied Physics Laboratory to solving the problems posed by V-SPELLS, our tool in development, called LUMOS, begins by applying methods from static analysis, natural language processing, and dynamic analysis to the legacy source code in order to generate high-level abstractions of these DSLs (Domain Specific Languages) that we call domain-specific semantic models (DSSMs) from the DSLs that comprise the source code.  These DSSMs will be generated in a language we refer to as the \textit{meta-DSL}, and in order to provide the patches to the legacy code requested in V-SPELLS, these DSSMs will have to be composed in very specific ways.  In order to ensure correctness of composition, as is required in V-SPELLS, we are providing verification via an algebraic framework using several ideas from category theory in which the key structure to our modeling is called an \textit{operad}.  Operads have begun to play an increasingly important role within applied mathematics (see \cite{LBPF,FBSD,BF,BO,GLMN}), and we find they provide an excellent mathematical model for our verification needs on V-SPELLS.  

To be more precise about our modeling, when a DSSM is written in the meta-DSL, we will use an operad to represent the DSSM in the meta-DSL, and composition of DSSMs within the meta-DSL will be modeled via a ``gluing" operation.  Mathematically, we are providing denotational semantics for a key portion of the language of the meta-DSL.  To ensure the highest level of correctness on composition between DSSMs in the meta-DSL, we aim provide a formalization of the denotational semantics of the meta-DSL.  In particular, we need to provide a formalization for the foundation for the denotational semantics of the meta-DSL:  operads.  We provide this formalization within the proof assistant Coq, and this is the focus of our paper.  

In Section \ref{infdefop}, we discuss the informal definition of operads; Section \ref{formodop} discusses the technicalities we faced and our solutions to defining operads within Coq; in Section \ref{forproofop}, we discuss our construction of the equivalent of an \textit{operad of sets} within Coq (namely, an \textit{operad of types}), and discuss our proof in Coq that this is an operad according to our specification in Section \ref{formodop}; and lastly, in Section \ref{relwork}, we compare our formalization to the only other formalization of operads we are aware of \cite{HH}.  We 

\section{Informally Defining Operads}\label{infdefop}
\blfootnote{\textbf{Distribution Statement A}:  Approved for Public Release, Distribution Unlimited}

While there does not seem to be an agreed-upon definition for a \textit{symmetric colored operad}, we note we are following the definition of a \textit{symmetric colored operad} in \cite{S1}.  However, we remark the definition in \cite{S1} does not include what is called the \textit{equivariance axiom} in \cite{MR3837179}; we too omit this axiom, since it is not relevant to what we want to accomplish in our work on V-SPELLS.  Regardless of these distinctions, we use \textit{operad} to mean \textit{symmetric colored operad} or \textit{colored operad} in the sequel.  

As our aim was to fully formalize the definition of an operad within Coq, we require precision, so we provide the full informal definition of an operad below in two parts.  The first part consists of the objects that comprise an operad.  

\begin{definition}[Data for an Operad]\label{defoperad1}
An \textbf{operad}, $\mco$, consists of a collection of types, which we will denote by $T$, and for each $n\geq 1$, $d\in T, \und c := c_0,\ldots, c_{n-1}$ a sequence of types in $T$, a collection of terms $\mco{d\choose \und c}$ such that, 
\begin{enumerate}
\item for each $c\in T$, we designate an element $\mbbo_c\in\mco{c\choose c}$ called the \textbf{$c$-colored unit};

\item if $\sigma$ is a permutation on $n$ letters, and $\und c \sigma:= c_{\sigma(0)},\ldots, c_{\sigma(n-1)}$, then there is a bijection between $\mco{d\choose \und c}$ and $\mco{d\choose \und c\sigma}$; 

\item for any sequence $\und b$ of types in $T$, if we denote by $\und c\bu_i \und b$ the sequence given by 

$$\underbrace{c_0,\ldots, c_{i-1}}_{\text{$\emptyset$ if $i = 0$}}, \und b, \underbrace{c_{i+1},\ldots, c_{n-1}}_{\text{$\emptyset$ if $i = n-1$}},$$

then we require the existence of a function:  

$$\circ_i: \mco{d\choose \underline c}\times \mco{c_i\choose\underline b} \rightarrow  
\mco{d\choose \und c\bu_i\und b}.$$

We typically refer to the function $\circ_i$ as \textbf{multi-composition}.  

\end{enumerate} 
\end{definition} 

\blfootnote{\textbf{Distribution Statement A}:  Approved for Public Release, Distribution Unlimited}
\begin{example}
For quick example of what the type signature of a multi-composition function looks like, let $\und c = c_0, c_1, c_2$, $\und b = b_0, b_1$, and $i = 1$, then $\circ_1$ has type signature:  

$$\mco{d\choose c_0, c_1, c_2}\times \mco{c_1\choose b_0, b_1}\ra \mco{d\choose c_0, b_0, b_1, c_2}.$$ 

\end{example} 

Now the data for an operad $\mco$ in Definition \ref{defoperad1} is subject to certain axiomatic constraints, and this forms the second half of our definition for an operad.  

\begin{definition}[Axioms for an Operad]\label{defoperad2} 

Let $\und c := c_0,\ldots, c_{n-1}, \und b := b_0,\ldots, b_{m-1}, \und a = a_0,\ldots, a_{\ell-1}$ be sequences from a collection of types $T$.  The axioms that the data for an operad $\mco$ must follow are given below.  

\begin{enumerate}
  \item The \textbf{horizontal associativity axiom}:  Suppose $n\geq 2$ and $0\leq i < j\leq n-1$, then for $(\alpha,\beta,\gamma)\in \mco{d\choose \und{c}}\times \mco{c_i\choose \und{a}}\times \mco{c_j\choose \und{b}}$,

$$(\alpha\circ_i \beta) \circ_{\ell-1+j} \gamma = (\alpha\circ_j \gamma) \circ_i \beta$$

To give a visual description of this axiom, we are requiring commutativity of the following diagram:  

\begin{diagram}
\mco{d\choose \und{c}}\times \mco{c_i\choose \und{a}}\times \mco{c_j\choose \und{b}} &\rTo^{(\circ_i, \text{id})} & \mco{d\choose \und{c}\bu_i \und a}\times \mco{c_j\choose \und{b}}\\
\dTo^{\cong}_{\text{swap}} &&  \dTo_{\circ_{\ell-1+j}}\\ 
\\
\mco{d\choose \und{c}} \times \mco{c_j\choose \und{b}}\times \mco{c_i\choose \und{a}}\\ 
\dTo^{(\circ_j, \text{id})} &&  \mco{d\choose (\und{c}\bu_i \und a)\bu_{\ell-1+j } \und b}\\
&& \parallel \\
\mco{d\choose \und{c} \bu_j \und b} \times \mco{c_i\choose \und a} & \rTo _{\circ_i} & \mco{d\choose (\und{c}\bu_j \und b)\bu_i \und a}. \\
\end{diagram}

\item The \textbf{vertical associativity axiom}:   Suppose $m, n\geq 1$, $0\leq i\leq n-1$, and $0\leq j\leq m-1$.  Then for $(\alpha,\beta,\gamma)\in \mco{d\choose \und{c}}\times \mco{c_i\choose \und{b}}\times \mco{b_j\choose \und{a}}$, 

$$(\alpha \circ_i \beta) \circ_{i+j} \gamma = \alpha\circ_i (\beta\circ_j \gamma)$$ 

That is, we are requiring commutativity of the following diagram:  

\begin{diagram}
\mco{d\choose \und{c}}\times \mco{c_i\choose \und{b}}\times \mco{b_j\choose \und{a}} & \rTo^{(\text{id}, \circ_j)} & \mco{d\choose \und c}\times \mco{c_i\choose \und b\bu_j \und a}\\
\\
\dTo^{(\circ_i, \text{id})} && \dTo_{\circ_i} \\
&& \mco{d\choose c \bu_i (\und b \bu_j \und a)}  \\
&& \parallel \\
\mco{d\choose \und c \bu_i \und b}\times \mco{b_j\choose \und a} & \rTo_{\circ_{i +j}} & \mco{d\choose (\und c \bu_i \und b) \bu_{i + j} \und a} 
\end{diagram}

\item The \textbf{left unity axiom} requires that for $\alpha\in\mco{d\choose \und c}$ with $n\geq 1$, $\mbbo_d \circ_1 \alpha = \alpha.$  
\blfootnote{\textbf{Distribution Statement A}:  Approved for Public Release, Distribution Unlimited}

\item The \textbf{right unity axiom} requires that for $n \geq 1$, $0 \leq i \leq n-1$, and $\alpha\in\mco{d\choose \und c}$, $\alpha \circ_i \mbbo_{c_i} = \alpha$.

\end{enumerate}
\end{definition} 

Before we give an example, some comments are in order about Definition \ref{defoperad2}.  

\begin{remark}\label{sancheck} 
\item We want to give some sanity checks of the associativity axioms.  First notice the following equality occurs in the right-hand corner of the diagram for the horizontal associativity axiom (1 of Definition \ref{defoperad2}):   

\begin{equation}\label{eq:op1}
\\\mco{d\choose (\und{c}\bullet_i \und a)\bullet_{\ell-1+j} \und b} = \mco{d\choose (\und{c}\bullet_j \und b)\bullet_i \und a}
\end{equation}

This equality arises from an equality of the following sequences:  

\begin{eqnarray} \label{eq:horiz}
(\und c \bullet_i \und a) \bullet_{\ell-1+j}\und b &=& (c_0,\ldots, c_{i-1}, \und a, c_{i+1},\ldots, c_{n-1})\bullet_{\ell-1+j} \und b\\ 
&=& c_0,\ldots, c_{i-1}, \und a, c_{i+1}, \ldots c_{j-1}, \und b, c_{j+1},\ldots, c_{n-1}\nonumber \\
&=& (\und c \bullet _j \und b) \bullet_i \und a \nonumber
\end{eqnarray}

In particular, in providing a specification in Coq for operads, we need to provide a proof that \eqref{eq:horiz} holds for such sequences in $T$.  

A similar equality of sequences is required to define the vertical associativity diagram:  

\begin{eqnarray} \label{eq:vert}
\und c \bullet_i (\und b \bullet_j \und a) &=& c_0,\ldots, c_{i-1}, (\und b \bullet_j \und a), c_{i+1},\ldots, c_{n-1}\\ 
&=& c_0,\ldots, c_{i-1}, b_0, \ldots, b_{j-1},\und a, b_{j+1},\ldots, b_{m-1} , c_{i+1},\ldots, c_{n-1}\nonumber\\
&=& (\und c \bullet _i \und b) \bullet_{i+j} \und a\nonumber
\end{eqnarray}  
\end{remark} 

While our definition seems extraordinarily abstract, the next example helps clarify the roots of the abstraction found in Definition \ref{defoperad1} and Definition \ref{defoperad2}.  Moreover, the next example will serve as the first application of our formal definition of operads, as we will prove in Coq that our realization of this example is an operad according to our specification.  

\begin{example}\label{exsets}
If we let $T$ be a collection of types for which $T$ is closed under finite products, we can define an operad $\sets_T$ by setting 

$$\sets_T{d\choose c_0,\ldots, c_{n-1}} := \hm(c_0\times\cdots\times c_{n-1}, d),$$
where the hom-set on the right is the collection of all functions from the product of sets $c_0\times\cdots\times c_{n-1}$ to the set $d$.  Given $c\in T$, the identity function on $c$ operates as the $c$-colored unit in $\sets_T{c\choose c} = \hm(c,c)$.  In this setting, we can explicitly define multi-composition $\circ_i$ from Definition \ref{defoperad1} which returns, given $f\in \hm(c_0\times\cdots\times c_{n-1}, d) $ and $g\in \hm(b_0\times\cdots\times b_{m-1}, c_i)$, the function $f\circ_i g$ which acts on the $(n+m-1)$-tuple $(x_0,\ldots, x_{i-1},\underline y,x_{i+1},\ldots, x_{n-1})$ as
$$(f\circ_i g)(x_0,\ldots, x_{i-1},\underline y,x_{i+1},\ldots, x_{n-1}) = f(x_0,\ldots, x_{i-1}, g(\underline y), x_{i+1},\ldots, x_{n-1}).$$ 

All other pieces of Definition \ref{defoperad1} and \ref{defoperad2} not mentioned above can be proved for $\sets_T$ using everything defined above and basic facts in set theory.  
\end{example}

\section{Formally Modeling Operads in Coq}\label{formodop}
\blfootnote{\textbf{Distribution Statement A}:  Approved for Public Release, Distribution Unlimited}
In defining the collection of terms $\mco{d\choose c_0,\ldots, c_{n-1}}$ in Coq, Definition \ref{defoperad1} requires that $d, c_i$ come from the collection $T$.  Throughout our specification in this paper, we will replace $T$ with one of Coq's in-house universes: $\type$.  In practice, we do need a proper subset of $\type$, but for simplicity in our paper, we use $\type$.  In the event we need a restriction to a subset of $\type$, we briefly discuss how to use Tarski universes to do this after the description of our formal model in Coq.  

\subsection{Encoding an Operad in Coq}

The first goal to tackle in defining an operad is giving a formal definition of $\mco{d\choose \und c}$.

\begin{note}[A Definition for $\mco{d\choose \und c}$ in Coq]\label{opdefcoq}
Informally, part of an operad $\mco$ is a collection of sets indexed by pair $d:\type$ and $\und c := c_0,\ldots, c_{n-1}:\,\textbf{list}\, \type$.  Since this is a \textit{collection of sets}, it would be natural to use a \textit{record} in Coq to make this definition.  To do so, we create a record in Coq, which we denote as $\textbf{Operad}$, whose single field is given by a function with type signature: $\type \ra \textbf{list}\,\type \ra \type$.  An instantiation of \textbf{Operad} will yield a function $\mco:  \textbf{list}\,\type \ra \type\ra \type$, so that $\mco{d\choose \und c}$ yields our desired collection of terms.  
\end{note} 

We give an example of our definition from Note \ref{opdefcoq}.  

\begin{example}\label{extype1} 
Our goal in Section \ref{forproofop} is to provide a version of $\sets_T$ in Example \ref{exsets} in Coq for $T = \type$; we will denote this operad by $\typebl$.  In Coq, if $\und c = c_0,\ldots, c_{n-1} : \textbf{list}\,\type$ and $d:\type$, then the following is definable in Coq via recursion:  

$$\typebl{d\choose \und c} := c_0\ra\cdots\ra c_{n-1}\ra d.$$

In particular, \textit{terms} of type $\typebl{d\choose \und c}$ are $n$-ary functions with codomain defined by $\und c$, and with return type $d$.    
\end{example} 

In the rest of our model in Coq, we also use a record to denote the data that comprises the operad (as in Definition \ref{defoperad1}) or the constraints the data is subject to (as in Definition \ref{defoperad2}).  Each piece in Definition \ref{defoperad1} and \ref{defoperad2} is a proposition that must be satisfied.  We first detail how the the data from Definition \ref{defoperad1} will be encoded as propositions within Coq.  

\blfootnote{\textbf{Distribution Statement A}:  Approved for Public Release, Distribution Unlimited}

\begin{note}[Data for an Operad in Coq]\label{operadnote1}
\begin{enumerate}
    \item the existence of a $c$-colored unit in $\mco$ ($1$ from Definition \ref{defoperad1}):  for all $c:\type$, there is a $\mbbo_c\in\mco{c\choose c}$; 
    
    \item the requirement that there is a bijection between $\mco{d\choose \und c}$ and $\mco{d\choose \und c\sigma}$ for a permutation $\sigma$ on $n$ letters ($2$ from Definition \ref{defoperad1}):  for all $d:\type$, $\und c, \und c' : \textbf{list}\, \type$ with the length $\und c$ at least $1$, and $\und c$ and $\und c'$ are permutations of one another, there is a bijection between $\mco{d\choose \und  c}$ and $\mco{d\choose \und c'}$; 
    
    \item the requirement for the existence of $\circ_i$ ($3$ from Definition \ref{defoperad1}); for all $i, n: \nats$, $d, c_i: \type$, $\und b, \und c: \textbf{list}\,\type$, if $\und c$ has length $n$, $1\leq n$, $i < n$, and the $n$th entry of $\und c$ is $c_i$, there is a function of type $\mco{d\choose \und c}\times \mco{c_i\choose \und b}\ra \mco{d\choose \und c\bu_i \und b}$.  
\end{enumerate}
\end{note} 

To make our implementation in Coq clear in Note \ref{operadnote1}, some remarks are in order about how to make the above precise within Coq:  

\begin{remark}\label{operadnote1remark}
\begin{enumerate}
    \item Any time \textit{bijection} is used in this context, we are referring to a bijection in $\type$.  That is, if $t, t': \type$, then there are functions $f: t\ra t'$, $f': t'\ra t$ such that $f\circ f' = \text{id}_{t'}$, and $f'\circ f = \text{id}_t$.  This is easily definable in Coq.  
    \item To create a proposition that two lists, $\und c, \und c'$, are permutations of one another in Coq, we can use Coq's built-in type \textbf{Permutation}.  This says that $\textbf{Permutation}\, \und c\, \und c': \prop$ (where $\prop$ is the type of all propositions in Coq).  
    \item The operation $\bullet_i$ on lists can be defined in Coq by taking the first $i$ entries of $\und c$, concatenating the list $\und b$, and then concatenating the the last $n-i-1$ entries of $\und c$ to the previous concatenation.  
    \item In $3$ of Note \ref{operadnote1}, we need the use of the $n$th function within Coq.  This function requires a default element as part of its arguments, which means we would need to choose a default element from $\type$ to use consistently throughout.  The choice we make in Coq is the \textit{unit} type, which is the type used to represent singleton sets.  
\end{enumerate}
\end{remark} 

We encode Definition \ref{defoperad2} into Coq in a similar manner using records, denoting this record by \textbf{OperadLaws}.  However, there is more caution to be had, due in part to the discussion in Remark \ref{sancheck}.  To demonstrate this caution, we discuss our modeling of of the horizontal associativity axiom within Coq in explicit detail below.    

\begin{note}[Axioms for an Operad in Coq]\label{operadnote2}
\begin{enumerate}
The horizontal associativity axiom in an operad (1 in Definition \ref{defoperad2}) can be defined in Coq by first listing a collection of parameters that we refer to as $P$:  

\begin{enumerate}
    \item[$\bu$] $n, m, \ell, i, j : \nats$; 
    \item[$\bu$] $d, c_i, c_j : \type$; 
    \item[$\bu$] $\und a, \und b, \und c : \textbf{list}\, \type$
    \item[$\bu$] $\alpha:\mco{d\choose\und c}, \beta: \mco{c_i\choose \und b}, \gamma: \mco{b_j\choose \und a}$
    \item[$\bu$] $2\leq n$, $1\leq m$, and $1\leq \ell$; 
    \item[$\bu$] $i < j$ and $j < n$; 
    \item[$\bu$] $\und c$ has length $n$, $\und b$ has length $m$, and $\und a$ has length $\ell$; 
    \item[$\bu$] the $i$th entry of $\und c$ is $c_i$ and the $j$th entry of $\und c$ is $c_j$; 
\end{enumerate}

Using what is now in $P$, we can give a proof that the $i$th entry of $\und c \bu_j \und b$ is $c_i$, and a proof that the $(\ell-1+j)$th entry of $\und c\bu_i \und a$ is $c_j$; we add these proofs to $P$.  With this update to $P$, we can state our formalization of the horizontal associativity axiom in Coq:  for all parameters that comprise $P$,  Equation \eqref{eq:horiz} in Remark \ref{sancheck} holds, and there exists a type casting function $\mcc_{\textbf{assoc}}$ such that 

$$\mcc_{\textbf{assoc}}\, P\, ((\alpha\circ_i \beta) \circ_{\ell-1+j}\gamma) = (\alpha\circ_j \gamma) \circ_i \beta$$
\end{enumerate} 
\end{note} 

The type-casting function $\mcc_{\textbf{assoc}}$ is necessary, since we have defined in Coq for each $d:\type$ and $\und c:\textbf{list}\,\type$, that $\mco{d\choose\und c}$be a type in $\type$, and the casting function provides a proof that the equality of types in Equation \eqref{eq:op1} holds.  However, the existence of $\mcc_{\textbf{assoc}}$ relies entirely on the proof of the equality of lists in Equation \eqref{eq:horiz}.  Now the equality in Equation \eqref{eq:horiz} requires a significant effort, and the most difficult part of formalizing this axiom is in providing its proof.  

Providing a formal specification of all other axioms in Definition \ref{defoperad2} to be inserted into the fields of of our record \textbf{OperadLaws} follows the same path as above:  

\begin{enumerate}
    \item carefully curate the correct collection $P$ of parameters needed for the axiom; 
    \item add in any proofs needed that can be deduced from everything currently in $P$;  
    \item show that any required equality of lists holds (this will be necessary for all axioms in Definition \ref{defoperad2}); 
    \item create the necessary casting function.  
\end{enumerate}

We have one last comment to make on the choices in our model.   

\blfootnote{\textbf{Distribution Statement A}:  Approved for Public Release, Distribution Unlimited}
\begin{remark}
In \cite{MR3837179} the definition for operads says that if $\und c = \emptyset$, the empty list of symbols coming from the collection $T$, then the symbol $\mco{d\choose \emptyset}$ still has meaning.  Notice in Definition \ref{defoperad1}, we do not allow the existence of of such a symbol since we require that the list $\und c$ is \textit{not} empty.  Our reason for doing so is that our main application relies on giving a version of Example \ref{exsets} in Coq.  Within $\sets_T$, if $\und c = \emptyset$, then the product of an empty list of sets is a singleton, $\stl \bu \str$, so that $\sets_T{d\choose \emptyset} \cong \sets_T{d\choose \stl \bu\str}$.  We can model this situation in Coq by letting $\und c$ be the list whose only entry is $\mbu : \type$, the unit type.    
\end{remark}

\subsection{Tarski Universes}

A solution to using a subset $T$ of $\type$ is to define $T$ in Coq as a \textit{Tarski universe}.  This defines $T : \type$, as well as an \textit{interpretation} that allows the terms of $T$ be regarded as \textit{codes} for actual types.  In this way, the type $T$ is a set together with an injective mapping to $\type$, which is exactly the data of a subset of $\type$.  Our approach to implementing this definition in Coq involves the following:  

\begin{enumerate}
    \item a type $\mcb$ in Coq with nullary constructors, we call the \textit{base types}, and whose terms we refer to as \textit{type sigils};
    \item the constructors that define the type $T$, which include:
        \begin{itemize}
            \item a constructor with signature $\textbf{Ty}:  \mcb \ra T$ which encodes the base types into $T$; 
            \item other constructors that may model products, such as $\textbf{p}:  T\ra T\ra T$, or $\textbf{fn}: T\ra T\ra T$, which can model functions;   
        \end{itemize}
    \item an assignment for $\mcb$ within $\type$, and a recursively-defined interpretation function $\el:  T\ra \type$ that assigns a value within $\type$ to each $t: T$.  
\end{enumerate}

We give an example of what this would look like explicitly.  

\begin{example}
We define our collection of base types $\mcb$ in Coq with the nullary constructors $N, U$, and $B$.  Within Coq, we create a function $\textbf{I}$ that interprets these type sigils:  $N$ is assigned to $\nats$, the type of natural numbers; $U$ to $\mbu$, the unit type; $B$ to $\mbb$, which is bool.    

If we want to model products and functions within in $T$, then we can define $T$ with constructors:  

\begin{enumerate}
    \item $\textbf{Ty}:  \mcb \ra T$; 
    \item $\textbf{p}:  T\ra T\ra T$; 
    \item $\textbf{fn}: T\ra T\ra T$.
\end{enumerate}

Now $\textbf{El}$ will provide the embedding into Coq via the following recursion:  

\begin{enumerate}
    \item $\textbf{El}\,(\textbf{Ty}\, t) \Rightarrow \textbf{I}\,t$ 
    \item $\textbf{El}\,(\textbf{p}\, t\, t') \Rightarrow \textbf{El}\, t \times \textbf{El}\, t'$ 
    \item $\textbf{El}\,(\textbf{fn}\, t\, t') \Rightarrow \textbf{El}\, t \ra \textbf{El}\, t'$ 
\end{enumerate}

For an explicit example of a code in $T$, we have $\textbf{p}\, (\textbf{Ty}\,N) (\textbf{Ty}\,N): T$, and via the embedding $\textbf{El}$, this is a model for $\nats\times\nats$ in Coq.

\end{example}  

\blfootnote{\textbf{Distribution Statement A}:  Approved for Public Release, Distribution Unlimited}
\section{A Proof Using Our Model}\label{forproofop}

Our goal in this section is to discuss the formal proof that the equivalent of Example \ref{exsets} in Coq, which we denote by $\typebl$ and define in Example \ref{extype1}, is an operad according to our model.  

To formally demonstrate that $\typebl$ is an operad, we first need a definition of the function in the only field of the record \textbf{Operad} (see Note \ref{opdefcoq}).  Next, in the record \textbf{OperadLaws}, we need to define, for $c:\type$, the $c$-colored units (1 of Definition \ref{defoperad1}), provide proofs that $\typebl{d\choose\und c}$ is invariant (up to injection) under reordering of $\und c$ (2 of Definition \ref{defoperad1}), define the multi-composition functions (3 of Definition \ref{defoperad1}), and show all axioms in Definition \ref{defoperad2} hold according.  Wrapping these assignments and proofs together provides a term of type \textbf{Operad} and \textbf{OperadLaws}, which gives our desired formal proof.  

In Example \ref{extype1}, we give the definition of the required function in \textbf{Operad} for $\typebl$:  given $\und c = c_0,\ldots, c_{n-1}:\textbf{list}\,\type$ and $d:\type$, we write:   

$$\typebl{d\choose \und c} := c_0\ra\cdots\ \ra c_{n-1}\ra d,$$

which is the type of $n$-ary functions with codomain defined by $\und c$ and return type $d$.  Our instantiation of \textbf{OperadLaws} for $\typebl$ will use this definition throughout.  

The right-hand side of $\typebl{d\choose \und c}$ is defined via a recursive function, which we denote as \textbf{arr} (short for \textit{arrow}), with type signature $\,\textbf{list}\, \type \ra \type \ra \type$.  In particular, we define $\textbf{arr}\, \emptyset \, d = d$ (where $\emptyset$ is the empty list).  

\subsection{Implementing the Data for Type in Coq} 

Next we discuss in a series of notes, the implementation of Definition \ref{defoperad1} for \textbf{Type} in Coq, as well as the tools that were developed for use in this implementation.   

\begin{note}[$c$-colored units in \textbf{Type}]
If $\und c$ has single entry $c:\type$, then $\typebl{c\choose \und c} = c\ra c$, which is the type of all functions with domain and range given by $c$.  Then $\mbbo_{c} :=\text{id}_c$, the identity function on $c$.  
\end{note} 

\begin{note}[$\typebl{d\choose \und c} \cong \typebl{d\choose \und c \sigma}$]
Our motivation is to provide the equivalent of Example \ref{exsets} within Coq, and the analogous isomorphism for the operad $\sets_T$ is, 

$$\hm(c_0\times\cdots\times c_{n-1}, d)\cong \hm(c_{\sigma(0)}\times\cdots\times c_{\sigma(n-1)}, d),$$

given $\und c = c_0,\ldots, c_{n-1}$ and a permutation $\sigma$ on $n$ letters.  The isomorphism in the context of Coq asks us to construct a bijection between the two sets above.  Following Definition \ref{operadnote1} and comments in Remark \ref{operadnote1remark}, we can translate this into Coq for $\typebl$ as:  for all $d:\type$, $\und c, \und c' : \textbf{list}\, \type$ with the length $\und c$ at least $1$, and $\textbf{Permutation}\,\und c\, \und c'$, there is a bijection between $\typebl{d\choose \und  c}, \typebl{d\choose \und c'} : \type$.  

We can prove this in Coq using induction on the length of $\und c$, along with some preceding lemmas about the behavior of function composition and isomorphism in this setting.  
\end{note} 

\begin{note}[$\circ_i$ for $\typebl$]\label{circi in Coq}
Lastly, we need to write (3) of Definition \ref{operadnote1} for $\typebl$ in Coq.  Writing $\typebl{d\choose \und c}$ in the curried form, as opposed to using a verbatim translation of $\sets_T{d\choose \und c}$ from Example \ref{exsets}, provides the needed flexibility, via partial application, to implement the multi-composition function $\circ_i$ for $\typebl$ in Coq with respect to (3) of Definition \ref{defoperad1}.  The most important piece of our definition in Coq  is that we define a recursive function $\textbf{compose}$ with type signature

$$ \textbf{arr}\, \und c' \, (t \ra t') \ra \textbf{arr}\, \und b\, t \ra \textbf{arr}\, \und c'\, (\textbf{arr}\, \und b\, t');$$

where $t, t' :\type$, and $\und c', \und b : \textbf{list}\, \type$.  If $\und c = c_0,\ldots, c_{i-1}, c_{i},c_{i+1},\ldots, c_{n-1}$, we let $\und c' = c_0,\ldots, c_{i-1}$, $t = c_i$, and $t' = c_i\ra c_{i+1} \ra\cdots \ra c_{n-1} \ra d$, and this gives us, provided the correct type inference is written in, the required multi-composition operator $\circ_i$ for $\typebl$.  
\end{note}

\blfootnote{\textbf{Distribution Statement A}:  Approved for Public Release, Distribution Unlimited}
\subsection{Implementing the Axioms for Type in Coq}

Definition \ref{defoperad1} provides the base to build the axioms of an operad, which are given in Definition \ref{defoperad2}, and we discuss our implementation in Coq of this here.  There are several axioms listed in Definition \ref{defoperad2}, and as in Note \ref{operadnote2}, we keep our discussion focused on the horizontal associativity axiom (1 of Definition \ref{defoperad2}), as the proof that $\typebl$ satisfies all other axioms in Definition \ref{defoperad2} follows from a similar procedure in Coq.  

Our first hurdle comes from noticing that our definition for $\circ_i$ in Note \ref{circi in Coq} is what we want mathematically, but that Coq does not automatically recognize the equality of types:  

$$\und c'\, (\textbf{arr}\, \und b\, t') = \textbf{arr} (\und c \bu_i\und b)\, d,$$

where $\und c = c_0,\ldots, c_{i-1}, c_{i},c_{i+1},\ldots, c_{n-1}$, $\und c' = c_0,\ldots, c_{i-1}$, $t = c_i$, $t' = c_i\ra c_{i+1} \ra\cdots \ra c_{n-1} \ra d$, so that  (see 3 of Definition \ref{defoperad1}) $\und c \bu_i b = c_0,\ldots, c_{i-1}, \und b,c_{i+1},\ldots, c_{n-1}$.  Since we require for $f:\typebl{d\choose \und c}, g:\typebl{c_i\choose \und b}$ that $f\circ_i g :\typebl{d\choose \und c \bu_i \und b}$, this presents an issue whose solution is, as in Note \ref{operadnote2}, a type casting function. 

The remainder of our formalization of $\typebl$ within Coq is a tour de force of type casting, and we discuss the tools we use in this proof below.  First, we give the formal definition we use for a type cast within Coq.  

\begin{definition}\label{typecast} 
Given $A, B: \type$, and an equation, $A = B$, a \textbf{type cast} $\mcc_{A = B}$ is a function such that for $a : A$, $\mcc_{A=B}\, a : B$.  
\end{definition}

\blfootnote{\textbf{Distribution Statement A}:  Approved for Public Release, Distribution Unlimited}

In order to manipulate the type casts that occur throughout our proof that $\typebl$ satisfies the axioms in Definition \ref{defoperad2}, we prove a handful of general facts about type casts in Coq, which we discuss below.  

\begin{note}[Type Casts for Definition \ref{defoperad2} in Coq]\label{tcnote} 
\begin{enumerate}
\item  The composition two type casts is a type cast: given equations of types $A = B$ and $B = C$, we have $\mcc_{B = C}\circ \mcc_{A = B} = \mcc_{A = C}$.
\item A type cast using an equation of types $A = A$ (i.e., a type cast between two types Coq recognizes as identical) is equal to the identity: $\mcc_{A = A}\, a = a$ for $a:A$.   
\item Two type casts between the same two types (i.e., both using equations of type $A = B$) are equal:  for all $a: A$, $\mcc_{A = B}\, a = \mcc'_{A = B}\, a$.  
\item Type casting a function and then applying it to an argument is the same as applying the original function to an argument that had been type cast:  if $f: B\ra C$ and $a: A$, then $ (\mcc_{A\ra B = B\ra C}\, f)\, a= f\, (\mcc_{A = B}\, a)$
\end{enumerate}
\end{note}

These facts smooth the process of showing $\typebl$ satisfies the operad axioms of Definition \ref{defoperad2}, as this involves manipulations of several type casts. For example, 2 and 3 of Note \ref{tcnote} ensure that it is not necessary to keep track of how these manipulations impact the \emph{equations} on which the type casts rely, since we need only that the \emph{types} involved match in order to show equality.

Now we discuss how to utilize the facts we demonstrated in Note \ref{tcnote}, by discussing their use in our proof the horizontal associativity axiom (1 of Definition \ref{defoperad2}) is satisfied in $\typebl$.  Our first step is to prove the following key lemma of equality involving the \textbf{compose} function:

\begin{equation}\label{eq:compHaa}
\textbf{compose}\, (\,\mcc\,(\,\textbf{compose}\, (\,\mcc\, \alpha\,)\, \beta\,)\,)\, \gamma =\\
\mcc\, (\,\textbf{compose}\, (\,\mcc\, (\,\textbf{compose}\, (\,\mcc\, \alpha\,)\, \gamma\,)\,)\,\beta \,).
\end{equation}

Here, $\alpha, \beta, \gamma$ are terms of the appropriate types in the operad $\typebl$, and the type casting equations have been suppressed from the notation.  We can show the equality in Equation \ref{eq:compHaa} by induction on the appropriate lists and several uses of 4 of Note \ref{tcnote}.  Notice that Equation \ref{eq:compHaa} is essentially the horizontal associativity axiom in $\typebl$, and this is the case:  to prove the horizontal associativity axiom, we manipulate type casts using our work in Note \ref{tcnote} until they match the equality in Equation \ref{eq:compHaa}.  

Demonstrating the remainder of the axioms in Definition \ref{defoperad2} for $\typebl$ within Coq follows a similar pattern:  show the pattern for the given axiom holds for \textbf{compose}, and then manipulate the type casts appropriately using Note \ref{tcnote} to arrive at the desired axiom.

\section{Related Work }\label{relwork}
In \cite{HH}, the authors present a formalization of a simpler type of an operad using Cubical Agda, which is an extension of Agda with Cubical Type Theory.  Cubical Type Theory is an alternative to Homotopy type theory that is more directly amenable to constructive interpretations, so fully understanding the implementation of operads in \cite{HH} requires a working knowledge of a variant of Homotopy type theory, as well as how to use its implementation in Agda.  

Moreover, Agda does not have significant automation, so showing, for example, our proof in Section \ref{forproofop} would require \textit{significantly} more work.  However, we do not think that our work would be impossible to translate into Agda, just require much more boilerplate code (e.g. handwriting structural induction tactics).  

We also want to compare what was formalized in our work to that of \cite{HH}.  What they use in \cite{HH} to refer to an operad is an operad with a collection of types $T$ for which $|T| = 1$.  In particular, $\mco{d\choose\und c}$ can be parametrized by the natural numbers, so we can write $\mcp(n) := \mco{d\choose \und c}$, if $|\und c| = n$, where $\und c = c_0,\ldots, c_{n-1}$, with $c_i = d$.  Moreover, there is a unique identity ($1\in\mcp(1)$), the functions $\circ_i$ have type $\mcp(n)\ra \mcp(m)\ra \mcp(n+m-1)$, and there is a significant simplification of the associativity axioms in Definition \ref{defoperad2}. We also note \cite{HH} they define their singly-colored operads to have the \textit{equivariance axiom} given in \cite{MR3837179}.  
\blfootnote{\textbf{Distribution Statement A}:  Approved for Public Release, Distribution Unlimited}
\bibliography{CALCO23Paper}

\end{document}